\documentclass{article}
\usepackage{amsfonts}

\usepackage[dvips]{graphicx}
\usepackage[english]{babel}
\usepackage{amsmath}
\usepackage{latexsym}
\usepackage{verbatim}
\usepackage{varioref}
\usepackage[latin1]{inputenc}
\usepackage[T1]{fontenc}
\usepackage{array}
\usepackage{delarray}
\usepackage{subfigure}
\usepackage{epsfig}

\newenvironment{proof}[1][Proof]{\noindent\textbf{#1 :} }{\ \rule{0.4em}{0.4em}}

\begin{document}

\title{Time-Changed Bessel Processes and Credit Risk\footnote{We thank Marc Yor for all the
important contributions to this paper. Any remaining errors are our own.}}
\author{Marc Atlan \\
Laboratoire de Probabilit\'{e}s\\
Universit\'{e} Pierre et Marie Curie \and Boris Leblanc \\
Equities $\&$ Derivatives Quantitative R$\&$D \\BNP Paribas}
\date{First Draft: December 16 2004; Last Draft: March 15 2006 }
\maketitle

\begin{abstract}

The Constant Elasticity of Variance (CEV) model is mathematically
presented and then used in a Credit-Equity hybrid framework. Next,
we propose extensions to the CEV model with default: firstly by
adding a stochastic volatility diffusion uncorrelated from the stock
price process, then by more generally time changing Bessel processes
and finally by correlating stochastic volatility moves to the stock
ones. Properties about strict local and true martingales in this study are discussed.
Analytical formulas are provided and Fourier and Laplace transform
techniques can then be used to compute option prices and
probabilities of default.

\end{abstract}

\section{Introduction}

\indent It has been widely recognized for at least a decade that the
option pricing theory of Black and Scholes (1973) and Merton (1973)
is not consistent with market option prices and underlying dynamics.
It has been noted that options with different strikes and maturities
have different implied volatilities. Indeed, markets take into
account in option prices the presence of skewness and kurtosis
 in the probability distributions of log returns. In
order to deal with those effects, one could use stochastic
volatility models (e.g. Heston (1993), Hull and White (1987) or
Scott (1987)). Another common alternative is to use a deterministic
time and stock price dependent volatility function, the so-called
local volatility to capture these effects. One would then build the
volatility surface by excerpting the values of this function from
option prices, thanks to the well-known Derman and Kani (1994) and
Dupire (1994) formula.

One of the first models developed after Black Merton Scholes (1973)
is the Constant Elasticity of Variance model pioneered by Cox (1975)
where the volatility is a deterministic function of the spot level;
This latter model is somehow an ancestor of local volatility models.
It has very interesting features since it suggests that common
stock returns are heteroscedastic and that volatilities implied by
the Black and Scholes formula are not constant, in other words skew
exists in this model. Another interesting property is that it takes
into account the so called "Leverage Effect" which considers the
effects of financial leverage on the variance of a stock: a stock
price increase reduces the debt-equity ratio of a firm and therefore
decreases the variance of the stock's returns (see for instance
Black (1976), Christie (1982) or Schwert (1989)). A last but not
least feature of this model is that it has a non-zero probability of
hitting $0$ and this could be of importance when one is interested
in modelling default by defining bankruptcy as the stock price
falling to $0$.

For the last few years, the credit derivatives market has become
more and more important and the issue of modeling default has grown,
giving birth to two main classes of models. The first class is the
structural models of the firm pioneered by Merton (1974) where
bankruptcy occurs if the asset value falls to a boundary determined
by outstanding liabilities. Other early work on such models was done
by Black and Cox (1976) and Geske (1977). The other class commonly
called reduced-form models is less ambitious than structural models.
They consider the time of default as an exogenous parameter that
they calibrate under a risk neutral probability to market data.
These models were developed by Artzner and Delbaen (1995), Jarrow
and Turnbull (1995), Duffie, Schroder and Skiadas (1996) and Madan
and Unal (1998).

The credit risk is also a component of the equity derivatives market
as it may appear in convertible bonds or more generally in Capital
Structure Arbitrage for people that embedded it from out-of-the
money puts. It is then clear that having a consistent modeling of
equity and credit is essential to eventually be able to manage those
cross-asset positions. Indeed, a market standard has been
developed during the last few years which involves a jump diffusion
dynamics for the stock price with a local probability of default for
the jump factor. This kind of model has been presented for instance
in Ayache, Forsyth and Vetzal (2003). An important drawback of this
modeling is that the stock has to jump to zero in order to default,
which isn't a realistic assumption as we can see on several
historical data and as argued in Atlan and Leblanc (2005).

The necessity to have stock price diffusions that don't jump to zero
in order to default and still have a non-zero probability of falling
to zero leads us to naturally consider CEV processes. Moreover, CEV
models have the advantage to provide closed-form formulas for
European vanilla options and for the probability of default. Those
computations were originally performed by Cox (1975) in the case
where the stock can default and by Emanuel and McBeth (1982) when
the stock never defaults. Then, one may want to add a stochastic
volatility process to the CEV diffusion in order to capture some
volatility features such as a smile or such as a more
realistic volatility term structure.
Finally, to get more dependency between the stock price
and the volatility, one may add some correlation.

Those guidelines lead us to study in section 2 the one-dimensional
marginals, the first-passage times below boundaries and the default
of martingality of Constant Elasticity of Variance processes, mainly
by relating those latest to Bessel processes. In section 3, we
propose a CEV model that is stopped at its default time and we
provide closed form formulas for European vanilla options, Credit
Default Swaps and Equity Default Swaps. Section 4 extends the
Constant Elasticity of Variance framework to a Constant Elasticity
of Stochastic Variance one by firstly adding a stochastic volatility
to the CEV diffusion and in a second time more generally consider
time-changed Bessel processes with a stochastic integrated time
change. Quasi-analytical formulas conditionally on the knowledge of
the law of the time change are provided for vanilla options and CDSs
and examples are given. Section 5 adds a correlation term to the
general time-changed power of Bessel process framework, once again
quasi analytical formulas conditionally on the knowledge of the
joint law of the time change and a process related to the rate of
time change are provided for probabilities of default and for
vanilla options, and computations for several examples are shown.
All the models proposed in this paper are true martingales and the
martingale property is carefully proven for the different
frameworks. Finally, section 6 concludes and presents
possible extensions of this work.\\

\noindent\textbf{Convention} \textit{For strictly negative
dimensions we define squared Bessel processes up to their first
hitting time of $0$ after which they remain at $0$.}\\

We set this convention because we wish to consider positive Bessel
processes. For a study of negative dimension Bessel processes with
negative values, we refer to G\"{o}ing-Jaeschke and Yor (2003).

\section{A Mathematical Study of CEV Processes}
\subsection{Space and Time Transformations}

A reason why Bessel processes play a large role in financial
mathematics is that they are closely related to widely used models
such as Cox, Ingersoll and Ross (1985), i.e. the CIR family of
diffusions for interest rates framework, such as the Heston (1993)
stochastic volatility model or even to the Constant Elasticity of
Variance model of Cox (1976). They are more generally related to
exponential of time-changed Brownian motions thanks to Lamperti
(1972) representations.\\

Let us now concentrate on the CIR family of diffusions: they solve
the following type of stochastic differential equations:
\begin{equation}
dX_t=(a-bX_t)dt+\sigma\sqrt{|X_t|}dW_t {\label{cir}}
\end{equation}
with $X_0=x_0>0$, $a\in\mathbb{R}_+$, $b\in\mathbb{R}$, $\sigma>0$
and $W_t$ a standard Brownian motion. This equation admits a strong
(e.g. adapted to the natural filtration of $W_t$) unique solution
that takes values in $\mathbb{R}_+$.\\

Let us remark that squared Bessel processes of dimension $\delta>0$
can be seen as a particular case of a CIR process with $a=\delta$,
$b=0$ and $\sigma=2$. We also recall that a Bessel process $R_t$
solves the following diffusion equation
\begin{equation*}
dR_t=dW_t+\frac{\delta-1}{2R_t}dt
\end{equation*}
where for $\delta=1$, the latter $\frac{\delta-1}{2R_t}dt$ must be replaced by a local time term.

One is now interested in the representation of a CIR process in
terms of a time-space transformation of a Bessel Process:
\newtheorem{BCIR}{Lemma}[section]
\begin{BCIR}
A CIR Process $X_t$ which solves equation (\ref{cir}) can be
represented in the following form:
\begin{equation}
X_t=e^{-bt} BESQ_{(\delta,x_0)}(\frac{\sigma^2}{4b}(e^{bt}-1))
\end{equation}
where $BESQ_{(\delta,x_0)}$ denotes a squared Bessel Process
starting from $x_0$ at time $t=0$ of dimension
$\delta=\frac{4a}{\sigma^2}$
\end{BCIR}

\begin{proof}
This lemma results from the identification of two continuous
functions $f$ and $g$ (with g strictly increasing and $g(0)=0$) such
as
\begin{equation*}
X_t=f(t)BESQ_{(\delta,x_0)}(g(t))
\end{equation*}
To do so, we apply It\^{o}'s formula and Dambis (1965),
Dubins-Schwarz (1965) theorem
\end{proof}\\

This relation is widely used in finance, for instance in Geman and
Yor (1993) or Delbaen and Shirakawa (2002).

Let us now introduce the commonly called CEV (Constant Elasticity of Variance), which
was introduced by Cox (1975, 1996) and that solves the following equation:
\begin{equation}
dX_t=\mu X_tdt+\sigma {X^{\alpha}_t} dW_t {\label{cev}}
\end{equation}
with $X_0=x_0>0$, $\alpha\in\mathbb{R}$, $\mu\in\mathbb{R}$,
$\sigma>0$ and $W_t$ a standard brownian motion.
\newtheorem{BCEV}[BCIR]{Lemma}
\begin{BCEV}
A CEV Process $X_t$ which solves equation (\ref{cev}) can be
represented as a power of a CIR process, indeed for
$\beta=2(\alpha-1)$, $1/X_t^{\beta}$ solves
\begin{equation}
d\big(\frac{1}{X^{\beta}_t}\big)=\big(a-b\frac{1}{X^{\beta}_t}\big)dt+\Sigma\sqrt{|\frac{1}{X^{\beta}_t}|}dW_t
\end{equation}
where $a=\frac{\beta(\beta+1)\sigma^2}{2}$, $b=\beta\mu$,
$\Sigma=-\beta\sigma$ and .
\end{BCEV}

\begin{proof}
This lemma is just an application of It\^{o}'s Lemma.
\end{proof}\\

As a consequence of Lemma 2.1 and Lemma 2.2, one obtains the
following representation for a CEV process:
\newtheorem{CEVr}[BCIR]{Proposition}
\begin{CEVr}
A CEV Process $X_t$ solution of equation (\ref{cev}) can be
represented in the following form:
\begin{equation}
X_t=e^{\mu t}
BESQ^{\frac{1}{2(1-\alpha)}}_{(\frac{2\alpha-1}{\alpha-1},x^{-2(\alpha-1)}_0)}\bigg(\frac{(\alpha-1)\sigma^2}{2\mu}(e^{2(\alpha-1)\mu
t}-1)\bigg)
\end{equation}
where $BESQ_{(\delta,x_0)}$ denotes a squared Bessel Process
starting from $x_0$ at time $t=0$ of dimension $\delta$.
\end{CEVr}

\subsection{Distributions and Boundaries}

We will now recall well known results about squared Bessel
processes and deduce some properties about CEV processes.

\subsection*{Path Properties}

\newtheorem{BesProp}[BCIR]{Proposition}
\begin{BesProp}
According to its dimension, the squared Bessel process has
different properties:\\
(i)   if $\delta\leq0$, {0} is an absorbing point.\\
(ii)  if $\delta<2$, \{0\} is reached a.s.\\
(iii) if $\delta\geq2$, \{0\} is polar.\\
(iv)  if $\delta\leq2$, BESQ is recurrent.\\
(v)   if $\delta\geq2$, BESQ is transient.\\
(vi)  if $0<\delta<2$, \{0\} is instantaneously reflecting.
\end{BesProp}

\begin{proof}
The proof can be found in Revuz and Yor (2001).
\end{proof}\\

As a consequence, one may give some properties of the CEV
diffusions. A topic of interest for the remaining of the paper is whether or not $\{0\}$ is
reached by a CEV process.

\newtheorem{CEVProp}[BCIR]{Proposition}
\begin{CEVProp}
According to the value of $\alpha$, the CEV diffusion has
different properties:\\
(i)   if $\alpha<1$, \{0\} is reached a.s.\\
(ii)  if $\alpha\leq\frac{1}{2}$, \{0\} is instantaneously reflecting.\\
(iii) if $\frac{1}{2}<\alpha<1$, \{0\} is an absorbing point.\\
(iv)  if $\alpha>1$, \{0\} is an unreachable boundary.
\end{CEVProp}

\begin{proof}
It is a consequence of the previous proposition and of Proposition
2.3.
\end{proof}\\

\subsection*{Distributional Properties}






It is important to notice that the law of a squared Bessel process
can be seen in terms of non-central chi-square density:
\newtheorem{Beschi}[BCIR]{Lemma}
\begin{Beschi}
For any $BESQ_{\delta,x}$, one has:
\begin{equation}
BESQ_{\delta,x}(t)\overset{(d)}{=}t V^{(\delta,\frac{x}{t})}
\end{equation}
where $V^{(a,b)}$ is a non-central chi-square r.v. with $a$ degrees
of freedom and non-centrality parameter $b\geq 0$. Its density is
given by:
\begin{equation}
f_{a,b}(v)=\frac{1}{2^{\frac{a}{2}}}\exp\bigg(-\frac{1}{2}(b+v)\bigg)v^{\frac{a}{2}-1}\sum_{n=0}^{\infty}\bigg(\frac{b}{4}\bigg)^n\frac{v^n}{n!\Gamma(\frac{a}{2}+n)}
\end{equation}
\end{Beschi}

\begin{proof}
This proof results from simple properties of Laplace transforms and
can be found for instance in Delbaen and Shirakawa (2002).
\end{proof}
\\

We leave to the reader the calculation of the CEV density in terms
of non-central chi-square distributions.

Let us recall a useful result for the remaining of the paper on the
moments of a squared Bessel process:
\newtheorem{Chimom}[BCIR]{Corollary}
\begin{Chimom}
If  $V^{(a,b)}$ is a non-central chi-square r.v. with $a$ degrees of
freedom and noncentrality parameter $b\geq 0$, then for any real
constants $c$ and $d$:
\begin{equation}
\mathbb{E}[(V^{(a,b)})^c\mathbf{1}_{\{V^{(a,b)}\geq
d\}}]=e^{-\frac{b}{2}}2^c\sum_{n\geq 0}\big(\frac{b}{2}\big)^n
\frac{\Gamma(n+\frac{a}{2}+c)}{n!\Gamma\big(\frac{a}{2}+n)}
G(n+\frac{a}{2}+c,\frac{d}{2}\big) {\label{espbes}}
\end{equation}
where $G$ is defined as follows:
\begin{equation*}
G(x,y)=\int_{z\geq
y}\frac{z^{x-1}e^{-z}}{\Gamma(x)}\mathbf{1}_{\{z>0\}}dz
\end{equation*}
\end{Chimom}

\begin{proof}
This calculation is a simple application of Lemma 2.6.
\end{proof}
\\

Finally, for the computations involved in this paper, one recalls
the two following identities on the complementary non-central
chi-square distribution function $Q$ that one can find in Johnson
and Kotz (1970):
\begin{eqnarray*}
Q(2z,2\nu,2\kappa)&=&\sum_{n\geq 1}g(n,\kappa)G(n+\nu-1,z)\\
1-Q(2\kappa,2\nu-2,2z)&=&\sum_{n\geq 1}g(n+\nu-1,\kappa)G(n,z)
\end{eqnarray*}
where $g(x,y)=-\frac{\partial G}{\partial y}(x,y)$.

\subsection*{First-Hitting Times}

We now concentrate on the first hitting time of $0$ by a Bessel
process. For this purpose, let us consider a Bessel Process $R$ of
index $\nu>0$ starting from $0$ at time $0$, then, one has:
\begin{equation}
L_1(R)\overset{(d)}{=}\frac{1}{2 Z_{\nu}}
\end{equation}
where $L_1(R)=\sup\{t>0, R_t=1\}$ and $Z_{\nu}$ is a gamma variable
with index $\nu$ that has the following density:
\begin{equation}
\mathbb{P}(Z_{\nu}\in
dt)=\frac{t^{\nu-1}e^{-t}}{\Gamma(\nu)}\mathbf{1}_{\{t>0\}} dt
\end{equation}
This result is due to Getoor (1979). Thanks to results on time
reversal (see Williams (1974), Pitman and Yor (1980) and Sharpe
(1980)), we have:
\begin{equation}
(\hat{R}_{T_0-u}; u<T_{0}(\hat{R}))\overset{(d)}{=}(R_u; u<L_1(R))
\end{equation}
where $\hat{R}$ is a Bessel Process, starting from $1$ at time $0$
of dimension $\delta=2(1-\nu)$ and $T_{0}(\hat{R})=\inf\{t>0,
\hat{R}_t=0\}$. As a consequence, one has:
\begin{equation}
T_0(\hat{R})\overset{(d)}{=}\frac{1}{2 Z_{\nu}}
\end{equation}
Using the scaling property of the Squared Bessel Process, one may
write:
\begin{equation}
T_0(BESQ^{\delta}_{x})\overset{(d)}{=}\frac{x}{2 Z_{\nu}}
{\label{deflt}}
\end{equation}
with $\delta=2(1-\nu)$.\\
\\
Hence, we are now able to state the proposition below:

\newtheorem{CEVdefault}[BCIR]{Proposition}
\begin{CEVdefault}
The probability of a CEV diffusion solution of equation (\ref{cev})
to reach $0$ at time $T$ with $\alpha<1$ is given by:
\begin{equation}
\mathbb{P}(T_0\leq T|X_0=x_0)=G(\frac{1}{2(1-\alpha)},\zeta_T){\label{probdefcev}}
\end{equation}
where $G$ and $\xi_T$ are defined as follows:
\begin{eqnarray}
G(x,y)&=&\int_{z\geq
y}\frac{z^{x-1}e^{-z}}{\Gamma(x)}\mathbf{1}_{\{z>0\}} dz\\
\zeta_T&=&\frac{\mu
x^{2(1-\alpha)}_0}{(1-\alpha)\sigma^2(1-e^{2(\alpha-1)\mu T})}
\end{eqnarray}
\end{CEVdefault}

\begin{proof}
This proof is just a consequence of Proposition 2.3 and equation
(\ref{deflt}).
\end{proof}

\newtheorem{CEVrem}[BCIR]{Remark}
\begin{CEVrem}
The calculation of the probability of default was originally done by
Cox (1975).
\end{CEVrem}


In order to compute first-passage times of scalar Markovian
diffusions below a fixed level, let us recall It\^{o} and McKean
(1974) results. If $(X_t,t\geq0)$ is scalar Markovian
time-homogeneous diffusion starting from $x_0$ at time $0$ of
infinitesimal generator $\mathcal{L}$ and that we define
$\tau_H=\inf\{t\geq0; X_t\leq H\}$ for $H<x_0$, then for any
$\lambda>0$, we have
\begin{equation*}
\mathbb{E}[e^{-\lambda\tau_H}]=\frac{\phi_{\lambda}(x_0)}{\phi_{\lambda}(H)}
\end{equation*}
where $\phi_{\lambda}$ is solution of the ODE
\begin{equation*}
\mathcal{L}\phi=\lambda\phi
\end{equation*}
with the following limit conditions:\\
$\lim_{x\rightarrow\infty}\phi_{\lambda}(x)=0$\\
$\textmd{If 0 is a reflecting boundary then
}\phi_{\lambda}(0+)<\infty$\\
$\textmd{If 0 is an absorbing boundary then
}\phi_{\lambda}(0+)=\infty$\\
\\

As a first example, let us now consider the first-hitting time below
a fixed level $0<y\leq x$ of a Bessel process $R_t$ of dimension
$\delta=2(\nu+1)$ starting from $x$ :
\begin{equation*}
\tau_y=\inf\{t\geq0; R_t\leq y \}
\end{equation*}
The law of $\tau_y$ (see It\^{o} and McKean (1974), Kent (1978) or
Pitman and Yor (1980)) is obtainable from the knowledge of its
Laplace transform $\mathcal{L}$. One has for any positive $\lambda$
\begin{eqnarray*}
\mathcal{L}(\lambda)&=&\mathbb{E}[e^{-\lambda\tau_y}]\\
&=&\frac{x^{-\nu}K_{\nu}(x\sqrt{2\lambda})}{y^{-\nu}K_{\nu}(y\sqrt{2\lambda})}
\end{eqnarray*}
where $\nu\in \mathbb{R}\backslash\mathbb{Z}$ and $K_{\nu}$ is a
Modified Bessel function defined as follows:
\begin{eqnarray*}
K_{\nu}(x)&=&\frac{\pi}{2\sin(\nu\pi)}(I_{-\nu}(x)-I_{\nu}(x))\\
I_{\nu}(x)&=&\sum_{k=0}^{\infty}\frac{(x/2)^{\nu+2k}}{k!\Gamma(\nu+k+1)}
\end{eqnarray*}\\

As a second example that will be useful for the computation of EDS
prices, let us write the infinitesimal generator of a CEV process:

\begin{equation*}
\mathcal{L}_{CEV}\phi=\sigma^2 x^{2\alpha}\frac{d^2\phi}{dx^2}+\mu
x\frac{d\phi}{dx}
\end{equation*}
that must solve
$$\mathcal{L}_{CEV}\phi=\lambda\phi$$
with the following conditions:\\
$\phi_{\lambda}$ is a decreasing function,
$\lim_{x\rightarrow\infty}\phi_{\lambda}(x)=0$\\
If $\alpha\leq\frac{1}{2}$, then $\phi_{\lambda}(0+)=\infty$\\
If $\frac{1}{2}<\alpha<1$, then $\phi_{\lambda}(0+)<\infty$\\
\\

We obtain the following result whose computations of the Laplace
transforms were originally performed by Davydov and Linetsky (2001):

\newtheorem{barriercev}[BCIR]{Proposition}
\begin{barriercev}
For a CEV process solution of (\ref{cev}) with $\alpha<1$ and
$\mu\neq0$, then
\begin{equation}
\phi_{\lambda}(x)=x^{\alpha-\frac{1}{2}}\exp\bigg(-\frac{\mu
x^{2(1-\alpha)}}{\sigma^2(1-\alpha)}\bigg)W_{k,m}\bigg(\frac{|\mu|
x^{2(1-\alpha)}}{\sigma^2(1-\alpha)}\bigg)
\end{equation}
where
\begin{equation*}
k=sgn(\mu)(\frac{1}{4(1-\alpha)}-\frac{1}{2})-\frac{\lambda}{2|\mu|(1-\alpha)}\indent \textmd{and}\indent m=\frac{1}{4(1-\alpha)}
\end{equation*}
and $W_{k,m}$ is a Whittaker function
\end{barriercev}

The definition of the Whittaker function can be found for instance
in Abramowitz and Stegun (1972).

\subsection{Loss of Martingality}

Let us now state a result on some martingale properties of Bessel
processes which play an essential role in pricing theory as is well known:

\newtheorem{martbes}[BCIR]{Theorem}
\begin{martbes}
Let $R_t$ be a Bessel process of dimension $\delta$ starting from
$a\neq0$, then:\\
(i) If\ $\delta\leq0$, $R_t^{2-\delta}$ is a true martingale up to the first hitting time of $0$.\\
(ii) If\ $0<\delta<2$, the process $R_t^{2-\delta}- L_t$ is a
martingale where $L_t$ is a continuous increasing process carried by
the zeros
of $(R_t, t\geq 0)$. \\
(iii) If\ $\delta=2$, $log(R_t)$ is a strict local martingale.\\
(iv) If\ $\delta>2$, $R_t^{2-\delta}$ is a strict local
martingale. Moreover, the default of martingality is
\begin{equation}
\gamma^{(\delta)}(t)=\mathbb{E}[R^{2-\delta}_0]-\mathbb{E}[R^{2-\delta}_t]=a^{2-\delta}\mathbb{P}^{(4-\delta)}_a(T_0\leq
t)
\end{equation}
where $\mathbb{P}^{\delta}_{a}$ is the law of
$(R^{(\delta)}_t,t\geq0)$.
\end{martbes}

\begin{proof}
(i) and (ii): Since $\{0\}$ is reached a.s., we need to apply
It\^{o}'s formula in a positive neighborhood of $0$. Let us consider
$\epsilon>0$. We have:
\begin{equation*}
(\epsilon+R^2_t)^{1-\frac{\delta}{2}}=(\epsilon+a^2)^{1-\frac{\delta}{2}}+(2-\delta)\int_0^t
(\epsilon+R^2_s)^{-\frac{\delta}{2}}R_sdW_s+\epsilon\delta(1-\frac{\delta}{2})\int_0^t
\frac{ds}{(\epsilon+R_s^2)^{\frac{\delta}{2}+1}}
\end{equation*}
Then, as $\epsilon$ tends to zero, it is easy to see the first
term of the right hand side is a true martingale for $\delta <2$ and
that the second term of the right hand side is increasing whose
support is the zeros of $(R_t, t\geq 0)$ when $\delta\geq0$.
If $\delta<0$, $T_t^{2-\delta}$ is a true martingale.\\
(iii): By applying It\^{o} formula, we obtain
\begin{equation*}
\log(R_t)=\log(R_0)+\int_0^t \frac{dW_s}{R_s}
\end{equation*}
We then see that $\log(R_t)$ is a local martingale. We prove that it
is a strict local martingale by first using the fact that
\begin{equation*}
\mathbb{P}^{\delta}_{a}=\mathbb{P}^{\delta}_{0}\ast\mathbb{P}^{0}_{a}
\end{equation*}
then writing that
\begin{eqnarray*}
\mathbb{E}[\log(R_t)]&=&\frac{1}{2}\mathbb{E}[\log(R_t^2)]=\frac{1}{2}\mathbb{E}[\log(BESQ_{2,0}(t)+BESQ_{0,a}(t))]\\
&\geq&\frac{1}{2}\mathbb{E}[\log(BESQ_{2,0}(t))]
\end{eqnarray*}
and finally since $BESQ_{2,0}(t)\overset{d}{=}2t\mathbf{e}$ where
$\mathbf{e}$ is a standard exponential, we obtain
\begin{equation*}
\mathbb{E}[\log(R_t)]\geq
C+\frac{1}{2}\log(t)\longrightarrow_{t\rightarrow\infty}+\infty
\end{equation*}
which shows that $\log(R_t)$ is not a true martingale.\\
(iv): To compute $\gamma^{(\delta)}$, we will need the following
result:
\newtheorem{girsbes}[BCIR]{Lemma}
\begin{girsbes}
Let $(R^{(\delta)}_t,t\geq0)$ be a Bessel process of dimension
$\delta>2$ starting from $a\neq0$, then
\begin{equation*}
\mathbb{P}^{4-\delta}_{a|\mathcal{R}_t\cap\{t<T_0\}}=\bigg(\frac{R^{(\delta)}_t}{a}\bigg)^{2-\delta}\cdot\mathbb{P}^{\delta}_{a|\mathcal{R}_t}
\end{equation*}
where $\mathcal{R}_t$ is the canonical filtration of the Bessel
process and $T_0$ the first-hitting time of the level $0$.
\end{girsbes}
\begin{proof}
This property results from a double application of Girsanov Theorem
by computing
\begin{equation*}
\frac{d\mathbb{P}^{4-\delta}_{a|\mathcal{R}_t\cap\{t<T_0\}}}{d\mathbb{P}^{2}_{a|\mathcal{R}_t}}\indent
\textmd{and}\indent
\frac{d\mathbb{P}^{\delta}_{a|\mathcal{R}_t}}{d\mathbb{P}^{2}_{a|\mathcal{R}_t}}
\end{equation*}
Then, by identification, one gets the announced result. A more
general result can be found in Yor (1992).
\end{proof}\\
\\
We may then write
\begin{equation*}
\mathbb{E}^{(\delta)}[R^{2-\delta}_t]=\mathbb{E}^{(4-\delta)}[a^{2-\delta}\mathbf{1}_{\{t<T_0\}}]
\end{equation*}
and consequently compute the default of martingality.
\end{proof}\\
\\
A proof in the case $0<\delta<2$ can be found in
Donati-Martin et al.\ (2006) and proofs when $\delta>2$ exist in
Elworthy, Li and Yor (1999). As a consequence, we obtain similar
results for a CEV process.

\newtheorem{martcev}[BCIR]{Proposition}
\begin{martcev}
Let $X_t$ be a CEV Process of elasticity $\alpha$ solving the following equation
\begin{equation*}
dX_t=\mu dt +\sigma X_t^{\alpha}dW_t
\end{equation*}
then:\\
(i) If\ $\alpha\leq\frac{1}{2}$, the process $e^{-\mu t}X_t$ is a
true martingale up to the first hitting time of $0$.\\
(ii) If\ $\frac{1}{2}<\alpha<1$, the process $e^{-\mu t}X_t- L^X_t$
is a martingale where $L^X_t$ is a continuous increasing process
carried by the zeros of $(X_t, t\geq 0)$ and consequently $e^{-\mu
t}X_t$ is a true martingale up to the first hitting time of $0$.\\
(ii) If\ $\alpha=1$, $e^{-\mu t}X_t$ is a geometric Brownian motion and hence a martingale.\\
(iii) If\ $\alpha>1$, $e^{-\mu t}X_t$ is a strictly local
martingale. Moreover, the default of martingality is
\begin{equation}
\gamma_{X}(t)=\mathbb{E}[X_0]-\mathbb{E}[e^{-\mu t}
X_t]=x_0G(\frac{1}{2(\alpha-1)},\zeta_T)
\end{equation}
where $G$ and $\zeta_T$ are defined as follows:
\begin{eqnarray*}
G(x,y)&=&\int_{z\geq
y}\frac{z^{x-1}e^{-z}}{\Gamma(x)}\mathbf{1}_{\{z>0\}} dz\\
\zeta_T&=&\frac{\mu
x^{2(1-\alpha)}_0}{(\alpha-1)\sigma^2(e^{2(\alpha-1)\mu T}-1)}
\end{eqnarray*}
\end{martcev}

\begin{proof}
This is just an application of Theorem 2.11, Proposition 2.3 and
equation (\ref{probdefcev}).
\end{proof}\\

\noindent A proof of the failure of the martingale property can be found in
Lewis (1998).

\newtheorem{pccev}[BCIR]{Remark}
\begin{pccev}
For $\alpha>1$, one has $\forall (t,K)\in\mathbb{R}^2_+$:
\begin{equation}
\mathbb{E}[(e^{-\mu t}X_t-K)_+]-\mathbb{E}[(K-e^{-\mu
t}X_t)_+]+\gamma_{X}(t)=\mathbb{E}[X_0]-K
\end{equation}
\end{pccev}

The last equation shows that in the case of a strictly local
martingale, the Call price must incorporate the default of
martingality in order to remain in a No Arbitrage model. For a study
on option pricing for strict local martingales, we refer to Madan
and Yor (2006) for continuous processes and to Chybiryakov (2006)
for jump-diffusion processes. Lewis (2000) also did this study in the case of explosions with stochastic volatility models
and in particular for a CEV diffusion.

\section{Credit-Equity Modelling}
\subsection{Model Implementation}

Usually, in the mathematical finance literature, one defines
a CEV diffusion for the stock price dynamics $S$ to be
\begin{equation*}
\frac{dS_t}{S_t}=\mu dt + \sigma S^{\alpha-1}_t dW_t
\end{equation*}
First of all, in a credit perspective, we will just consider the
case $\alpha<1$ since we are interested in models with a non-zero
probability of default. Once the stock has reached zero, the firm
has bankrupted and that is the reason why we stop the CEV diffusion
at its first default time. Then from what has been proven above, we
know that the stock price process hence defined is a true martingale
and that ensures the Absence of Arbitrage and moreover the uniqueness of the solution.
Hence, the stock price diffusion now becomes under the risk-neutral pricing measure:
\begin{equation*}
\frac{dS_t}{S_t}=r dt + \sigma S^{\alpha-1}_t dW_t\indent
\textmd{if}\indent t<\tau.
\end{equation*}
\begin{equation*}
S_t=0\indent \textmd{if}\indent t\geq\tau.
\end{equation*}
where $\tau=T_0(S)=\inf\{t>0, S_t=0\}$. In other words, the stock
price process considered is nothing else than a stopped CEV
diffusion $(S_{t\wedge\tau})_{t\geq 0}$.

\newtheorem{AOArem}[BCIR]{Remark}
\begin{AOArem}
Delbaen and Shirakawa (2002) showed the existence of a risk-neutral
probability measure whose uniqueness is only ensured on the stock
price filtration considered at time $\tau$
$\mathcal{F}_{\tau}=\sigma(S_t, t\leq\tau)$. Since our purpose is to
compute the price of options whose payoffs are
$\mathcal{F}_{\tau}-measurable$, we have the uniqueness of the
no-arbitrage probability.
\end{AOArem}

\subsection{European Vanilla Option Pricing}

Lemma 2.12 states that
\begin{equation}
\mathbb{P}^{4-\delta}_{x|\mathcal{R}_t\cap\{t<T_0\}}=\bigg(\frac{R^{(\delta)}_t}{x}\bigg)^{2-\delta}\cdot\mathbb{P}^{\delta}_{x|\mathcal{R}_t}
{\label{probid}}
\end{equation}
Thanks to this identity, we obtain the law of the stopped CEV
diffusion at a given time. Lemma 2.6 and Corollary 2.7 enable us to
compute the call and put option price:

For the call $C_0$ option price
\begin{eqnarray*}
C_0&=&e^{-rT}\mathbb{E}[(S_{T\wedge
\tau}-K)_+]\\&=&e^{-rT}\mathbb{E}[(S_T-K)_{+}\mathbf{1}_{T<\tau}]
\end{eqnarray*}and the put $P_0$ option price:
\begin{eqnarray*}
P_0&=&e^{-rT}\mathbb{E}[(K-S_{T\wedge\tau})_{+}]\\&=&e^{-rT}\mathbb{E}[(K-S_T)_{+}\mathbf{1}_{T<\tau}]+
K e^{-rT}\mathbb{P}(\tau\leq T)
\end{eqnarray*}Consequently, for the call price:
\begin{equation*}
C_0=S_0
Q(z_T,2+\frac{1}{1-\alpha},2\zeta_T)-Ke^{-rT}(1-Q(2\zeta_T,\frac{1}{1-\alpha},z_T))
\end{equation*}
and for the put price:
\begin{eqnarray*}
P_0&=&Ke^{-rT}\big(Q(2\zeta_T,\frac{1}{1-\alpha},z_T)-G(\frac{1}{2(1-\alpha)},\zeta_T)\big)\\&&-S_0
(1-Q(z_T,2+\frac{1}{1-\alpha},2\zeta_T))+K
e^{-rT}\mathbb{P}(\tau\leq
T)\\
&=&Ke^{-rT}Q(2\zeta_T,\frac{1}{1-\alpha},z_T)-S_0
(1-Q(z_T,2+\frac{1}{1-\alpha},2\zeta_T))
\end{eqnarray*}
where
\begin{eqnarray*}
z_T&=&\frac{2 r
K^{2(1-\alpha)}}{\sigma^2(1-\alpha)(e^{2(1-\alpha)r T}-1)}\\
\zeta_T&=&\frac{r
S_0^{2(1-\alpha)}}{(1-\alpha)\sigma^2(1-e^{-2(1-\alpha)r T})}
\end{eqnarray*}

Hence, one easily verifies that the put-call parity is satisfied.
Closed-form CEV  option pricing formulas were originally computed by
Cox (1975) for $\alpha<1$ and Schroder (1989) expressed those
formulas in terms of non-central chi-square distributions. Computing
option prices using the squared Bessel processes distributions was
done by Delbaen and Shirakawa (2002).

\subsection{Pricing of Credit and Equity Default Swaps}

Since we are dealing with default probabilities, it is obvious to
consider derivative products relying on these probabilities. One of
the most liquid protection instruments against default is the Credit
Default Swap (CDS). The buyer of the protection agrees to pay
periodical amounts until a default time (if it occurs) and in
exchange receives a cash amount which is a notional amount minus a
recovery rate in the case the company on which the contract is
written, defaults. The payoff of such kind of contract is:
\begin{equation*}
\Pi_{CDS}=-\sum_{i=1}^{n}e^{-rT_i}C\mathbf{1}_{\{\tau>T_i\}}+e^{-r\tau}(1-R)\mathbf{1}_{\{\tau\leq
T_n\}}
\end{equation*}
where C is the periodical coupon, $T_1,...,T_n$ the payment dates, R
the recovery rate assumed to be deterministic and $\tau$ the default
time. For simplicity purposes, we consider in this paper
deterministic interest rates. The CDS Fair Price is the expectation
of the payoff conditionally to the spot price filtration taken at
the pricing time, e.g.:
\begin{equation*}
CDS_t(T_1,T_n;C;R)=-C\sum_{i=1}^{n}e^{-r(T_i-t)}\mathbb{P}(\tau>T_i|S_t)+(1-R)\mathbb{E}[e^{-r(\tau-t)
}\mathbf{1}_{\{\tau\leq T_n\}}|S_t]
\end{equation*}
By absence of arbitrage, one must have $CDS_t(T_1,T_n;C;R)=0$ and
then
\begin{equation*}
C=\frac{(1-R)\mathbb{E}[e^{-r(\tau-t)}\mathbf{1}_{\{\tau\leq
T_n\}}|S_t]}{\sum_{i=1}^{n}e^{-r(T_i-t)}\mathbb{P}(\tau>T_i|S_t)}
\end{equation*}
From Proposition 2.8, we know the value of
$(\mathbb{P}(\tau>T_i|S_t)_{1\leq i\leq n})$. It then remains to
compute the following quantity
$\mathbb{E}[e^{-r\tau}\mathbf{1}_{\tau\leq t}]$ to be able to price
the CDS coupon C. By an integration by parts, we show that
\begin{equation}
\mathbb{E}[e^{-r\tau}\mathbf{1}_{\tau\leq
t}]=e^{-rt}\mathbb{P}(\tau\leq t)+r\int_0^t e^{-rs}\mathbb{P}(\tau
\leq s)ds
\end{equation}
Otherwise, one could just obtain this expectation by directly using the density of the first-hitting
time of $0$ that is provided  by the differentiation of the cumulative distribution function
:
\begin{equation*}
f_{\tau}(t)=\frac{2r(1-\alpha)\zeta_t^{\frac{1}{2(1-\alpha)}}e^{-\zeta_t}}{\Gamma(\frac{1}{2(1-\alpha)})(e^{2(1-\alpha)rt}-1)}
\end{equation*}
where $\zeta_t$ is defined above.\\
EDSs are very similar to CDSs except that payouts occur when the
stock price falls under a pre-defined level, which is often referred
to as a trigger price. The trigger price is generally between 30$\%$
and 50$\%$ of the equity stock price at the beginning of the
contract. Hence, these contracts provide a protection against a
credit event happening on the equity market for the buyer. They were
initiated by the end of 2003. At that time, it had become difficult
in many countries to structure investment-grade credit portfolios
with good returns because the CDS spreads were tightening, as
reported by Sawyer (2003). Another reason why people have
interest in those contracts is because the settlement
of the default is directly observed on the stock price.
Let us now define $\tau_{L}$ as the first
passage time of the stock price process under the level $L <
S_0$. Formally, we write $\tau_{L} = \inf\{t > 0; S_t \leq L\}$. We
recall the general valuation formula of an EDS:
\begin{equation*}
EDS_t(T_1,T_n;C;R)=-C\sum_{i=1}^{n}e^{-r(T_i-t)}\mathbb{P}(\tau_L>T_i|S_t)+\mathbb{E}[e^{-r(\tau_L-t)
}\mathbf{1}_{\{\tau_L\leq T_n\}}|S_t]
\end{equation*}
where C is the coupon, $T_1,...,T_n$ the payment dates and $r$ the
risk-free interest rate. Again, by absence of arbitrage, we can find
the coupon price, by stating that at the initiation of the contract:
\begin{equation*}
EDS_{t=0}(T_1,T_n;C;R)=0
\end{equation*}
Or equivalently
\begin{equation*}
C=\frac{\mathbb{E}[e^{-r(\tau_L-t)}\mathbf{1}_{\{\tau_L\leq
T_n\}}|S_t]}{\sum_{i=1}^{n}e^{-r(T_i-t)}\mathbb{P}(\tau_L>T_i|S_t)}
\end{equation*}
In order to price the coupon $C$, one needs to evaluate:
\begin{equation*}
\mathbb{E}[e^{-r\tau_L}\mathbf{1}_{\{\tau_L\leq t\}}]\indent
\textmd{and}\indent\mathbb{P}(\tau_L\leq t)
\end{equation*}
An integration by parts gives the Laplace transform of
$\mathbb{P}(\tau_L\leq t)$ for any $\lambda>0$
\begin{equation*}
\int_0^{+\infty}dt e^{-\lambda t}\mathbb{P}(\tau_L\leq
t)=\frac{\mathbb{E}[e^{-\lambda\tau_L}]}{\lambda}
\end{equation*}
Applying Fubini theorem, one observes that
\begin{equation*}
\int_0^{+\infty}dt e^{-\lambda
t}\mathbb{E}[e^{-r\tau_L}\mathbf{1}_{\{\tau_L\leq
t\}}]=\frac{\mathbb{E}[e^{-(r+\lambda)\tau_L}]}{\lambda}
\end{equation*}

Hence using Proposition 2.10, one is able to compute the Laplace
transform of the desired quantities necessary to evaluate an EDS.
One can then use numerical techniques (see Abate and Whitt (1995)
for instance) to inverse the Laplace transform in order to evaluate
prices.

\section{Stochastic Volatility for CEV Processes}
\subsection{A Zero Correlation Pricing Framework}

\subsection*{Impact of a Stochastic Time Change}

Due to the very important dependency between the probability of
default, the level of volatility and the skewness, we were naturally
brought to consider extensions of the CEV model that could relax the
high correlation between these three effects. More precisely, in a
CEV model, if one first calibrates the implied at-the-money
volatility, then either the skewness or the CDS will be calibrated
on adjusting the elasticity parameter. Hence, to be able to get some
freedom on the volatility surface, a possible extension is to
introduce a stochastic volatility in the CEV model instead of a
constant volatility. A CEV diffusion with a stochastic volatility is
actually just a power of a squared Bessel Process with a stochastic
time change instead of having a deterministic one like in
Proposition 2.3.

Another extension is to consider a power of a Bessel Process time
changed by an independent increasing process. More precisely, one
writes the following process for the stock price:
\begin{equation}
S_t=e^{rt}BESQ^{1-\frac{\delta}{2}}_{(\delta,x)}(\xi_t)\indent
if\indent t<\tau{\label{bestceq}}.
\end{equation}
\begin{equation*}
S_t=0\indent if\indent t\geq\tau.
\end{equation*}
where $x=S_0^{\frac{2}{2-\delta}}$, $\tau=T_0(S)=\inf\{t>0,
S_t=0\}=\xi^{-1}(T_0(BESQ))$ and $\xi_t$ is an strictly increasing
continuous integrable process independent from the squared Bessel
process. Subordinating a continuous process by an independent
L\'{e}vy process is an idea that goes back to Clark (1973).
Stochastic time changes are somehow equivalent to adding a
stochastic volatility in stock price diffusions. The basic intuition
underlying this approach could be foreseen through the scaling
property of the Brownian motion, or through Dambis (1965),Dubins and
Schwarz (1965) (DDS) theorem or even its extension to
semimartingales by Monroe (1978). More recently, Carr et al. (2003)
generated uncertainty by speeding up or slowing down the rate at
which time passes with a L\'{e}vy process. Our approach differs from
the one done in the L\'{e}vy processes literature for
mathematical finance: We are not considering the
exponential of a time changed L\'{e}vy process but a power of a time
changed Bessel process. Thanks to Lamperti representation (1972),
this means that we are considering a time changed geometric Brownian
motion $B$. More precisely, it is known that
\begin{equation*}
R_t=\exp\big(B_{C_t}+\nu C_t\big)\indent \textmd{and} \indent
C_t=\int_0^t\frac{ds}{R^2_s}
\end{equation*}
where $(R_t,t\geq 0)$ is a Bessel process of dimension $\delta=2(1+\nu)$ starting from $a\neq 0$.
Hence the time change considered in the stock price is
\begin{equation*}
Y_t=\int_0^{\xi_t}\frac{ds}{R^2_s}
\end{equation*}
and the stock price process as defined in equation (\ref{bestceq}) can be identified as follows:
\begin{equation*}
S_t=e^{rt}\exp\big(-2\nu B_{Y_t}-\frac{(2\nu)^2}{2}Y_t\big)
\end{equation*}
As a consequence, we have now proposed a new class of time changes where analytical computations are possible thanks to
a good knowledge of Bessel processes.\\
For the absence of arbitrage property, there must exist
a probability under which all the actualized stock prices are martingales.
A very simple property on martingales is that a process $M_t$ is a martingale if and only
if for every bounded stopping time $T$,
$\mathbb{E}[M_T]=\mathbb{E}[M_0]$. Nonetheless, this result is not
very convenient. \noindent Let us state and give a straightforward proof of the
martingality of the stock price process

\newtheorem{btcmart}[BCIR]{Proposition}
\begin{btcmart}
Consider
$M_t=BESQ^{1-\frac{\delta}{2}}_{(\delta,x)}(\xi_{t\wedge\tau})$
where following the previous hypotheses $\xi_t$ is a strictly
increasing continuous integrable process independent from $BESQ$,
$\tau$ is the $(M_t, t\geq0)$ first hitting time of $0$ and
$BESQ_{(\delta,x)}$ is a squared Bessel process of dimension
$\delta$ starting from $x\neq0$, then $(M_t, t\geq0)$ is a true
martingale.
\end{btcmart}

\begin{proof}
Let us define $\mathcal{R}_t=\sigma(R_s;s\leq t)$. We then naturally write the canonical
filtrations
$\mathcal{R}_{\xi_t}=\sigma(R_{\xi_s};s\leq t)$ and $\Xi_t=\sigma(\xi_s;s\leq t)$. For any bounded functional $F$,
we want to compute
\begin{equation*}
\mathbb{E}[F\big(R_{\xi_u};u\leq s\big)\big(R_{\xi_t}^{2-\delta}-R_{\xi_s}^{2-\delta}\big)]
\end{equation*}
Since $\xi$ is integrable and independent from $R$, we obtain by using Fubini theorem
\begin{eqnarray*}
\mathbb{E}[F\big(R_{\xi_u};u\leq s\big)\big(R_{\xi_t}^{2-\delta}-R_{\xi_s}^{2-\delta}\big)]
&=&\mathbb{E}\bigg[\mathbb{E}[F\big(R_{\xi_u};u\leq s\big)\big(R_{\xi_t}^{2-\delta}-R_{\xi_s}^{2-\delta}\big)\big|~\Xi_t]\bigg]\\
&=&\int\mathbb{P}_{\Xi_t}(da)\mathbb{E}[F\big(R_{a(u)};u\leq
s\big)\big(R_{a(t)}^{2-\delta}-R_{a(s)}^{2-\delta}\big)]
\end{eqnarray*}
The latest quantity is null by Theorem 2.11 and we have then shown
that for $s\leq t<\tau$
\begin{equation*}
R_{\xi_s}^{2-\delta}=\mathbb{E}[R_{\xi_t}^{2-\delta}|\mathcal{R}_{\xi_s}]
\end{equation*}
which is the announced result.
\end{proof}

\subsection*{Pricing Vanilla Options} One can find closed-form formulas for the call and put options
prices. Let us define the two following quantities $C_0(x,\delta,K,T;S_0)$
and $P_0(x,\delta,K,T;S_0)$:
\begin{eqnarray*}
C_0(x,\delta,K,T;S_0)&=&S_0
Q(\frac{(Ke^{-rT})^{\frac{2}{2-\delta}}}{x},4-\delta,\frac{S_0^{\frac{2}{2-\delta}}}{x})\\&&-Ke^{-rT}(1-Q(\frac{S_0^{\frac{2}{2-\delta}}}{x},2-\delta,\frac{(Ke^{-rT})^{\frac{2}{2-\delta}}}{x})) \\
P_0(x,\delta,K,T;S_0)&=&Ke^{-rT}Q(\frac{S_0^{\frac{2}{2-\delta}}}{x},2-\delta,\frac{(Ke^{-rT})^{\frac{2}{2-\delta}}}{x})\\&&-S_0
(1-Q(\frac{(Ke^{-rT})^{\frac{2}{2-\delta}}}{x},4-\delta,\frac{S_0^{\frac{2}{2-\delta}}}{x}))
\end{eqnarray*}
From there, one may obtain the option prices under the new general
framework.

\newtheorem{optprice}[BCIR]{Proposition}
\begin{optprice}
If one has the following stock price process:
\begin{equation*}
S_t=e^{rt}BESQ^{1-\frac{\delta}{2}}_{(\delta,x)}(\xi_t)\indent
if\indent t<\tau.
\end{equation*}
\begin{equation*}
S_t=0\indent if\indent t\geq\tau.
\end{equation*}
where $x=S_0^{\frac{2}{2-\delta}}$, $\tau=T_0(S)=\inf\{t>0, S_t=0\}$
and $\xi_t$ is a strictly increasing continuous integrable process
independent from $BESQ$ whose probability measure is
$\mu_{\xi_t}(dx)$, then:
\begin{eqnarray*}
C_0&=&\int_{\mathbb{R}_+}C_0(x,\delta,K,T;S_0)\mu_{\xi_T}(dx)\\
P_0&=&\int_{\mathbb{R}_+}P_0(x,\delta,K,T;S_0)\mu_{\xi_T}(dx)
\end{eqnarray*}
\end{optprice}

\begin{proof}
Let us prove this result for the call option price, a similar result
may be obtained for the put price. One has:
\begin{eqnarray*}
C_0&=&e^{-rT}\mathbb{E}[(S_T-K)_+]\\&=&e^{-rT}\mathbb{E}\big(\mathbb{E}[(S_T-K)_+|\sigma(\xi_s;
s\leq T)]\big)\\&=&\mathbb{E}[C_0(\xi_T,\delta,K,T;S_0)]
\end{eqnarray*}
\end{proof}

\subsection*{Computing the Default} Having the integrability of the
change of time and knowing its density, one could find a closed-form
formula for the probability of default $\tau=T_0(S)=\inf\{t>0,
S_t=0\}$ where
$S_t=e^{rt}BESQ^{1-\frac{\delta}{2}}_{(\delta,x)}(\xi_t)$. Let us
now compute the probability of default the proof of which is left to
the reader:
\newtheorem{optprob}[BCIR]{Proposition}
\begin{optprob}
If one considers a stock price process defined as follows:
\begin{equation*}
S_t=e^{rt}BESQ^{1-\frac{\delta}{2}}_{(\delta,x)}(\xi_t)
\end{equation*}
then the probability of default $\tau=\inf\{t>0, S_t=0\}$ is given by
\begin{equation*}
\mathbb{P}(\tau\leq T)=\mathbb{E}\bigg[G(1-\frac{\delta}{2},\frac{S_0^{\frac{2}{2-\delta}}}{\xi_T})\bigg]
\end{equation*}
where G is the complementary Gamma function.
\end{optprob}

\subsection{CESV Models}

Stochastic volatility models were used in a Black and Scholes (1973)
and Merton (1973) framework mainly to capture skewness and kurtosis
effects, or in terms of implied volatility skew and smile. In a
Constant Elasticity of Variance framework, one would use stochastic
volatility not to capture the leverage effect which partly already
exists due to the elasticity parameter but to obtain environments
for instance of low volatilities, high probabilities of default and
low skew. Let us consider an integrable jump-diffusion process
$(\sigma_t,t\geq0)$ to model the volatility. We will call those
diffusions Constant Elasticity of Stochastic Variance (CESV) for the
remainder of the paper. Leblanc (1997) introduced stochastic
volatility for CEV processes.

Hence, the class of models under a risk-neutral probability measure
proposed is of the following form:
\begin{equation*}
\frac{dS_t}{S_t}=r dt + \sigma_{t-} S^{\alpha-1}_t dW_t
\end{equation*}
where $\sigma$ is assumed to be independent from the Brownian motion driving the stock price
returns. Next, within an equity subject to bankruptcy framework, we
are going to stop the diffusion when the stock reaches $0$ just as
in the previous section. As a consequence, our diffusion becomes:
\begin{eqnarray*}
\frac{dS_t}{S_t}&=&r dt + \sigma_{t-} S^{\alpha-1}_t
dW_t\indent
if\indent t<\tau.\\
S_t&=&0\indent if\indent t\geq\tau.
\end{eqnarray*}
where $\tau=T_0(S)=\inf\{t>0, S_t=0\}$.\\
\noindent Before giving any concrete examples, let us show how CESV
models can be seen as Bessel processes with a stochastic time
change. So first, let us recall elementary results:

\newtheorem{tchange}[BCIR]{Lemma}
\begin{tchange}
Let R be a time change with $s\mapsto R_s$ continuous, strictly
increasing, $R_0=0$ and $R_t<\infty$, for each $t\geq0$, then for
any continuous semimartingale $X$ and any caglad (left continuous
with right limits) bounded adapted process $H$, one has:
\begin{equation}
\int_0^{R_t} H_s dX_s=\int_0^t H_{R_u} dX_{R_u}
\end{equation}
\end{tchange}

\begin{proof}
The proof can be found in Revuz and Yor (2001).
\end{proof}\\
\\
Then, using Lemma 4.4, (DDS) theorem and It\^{o} formula, we obtain
that
\begin{eqnarray*}
S_t&\overset{d}=&e^{rt}BESQ^{\frac{1}{2(1-\alpha)}}_{(2-1/(1-\alpha),S_0^{2(1-\alpha)})}(H_{t\wedge\tau})\\
\tau&=&\inf{\{t\geq0,S_t=0\}}\\
H_t&=&(1-\alpha)^2\int_{0}^{t}\sigma^2_s e^{-2(1-\alpha) r s}ds
\end{eqnarray*}
$H_t$ is by construction an increasing continuous integrable
process.\\ Hence $(e^{-rt}S_t,t\geq0)$ is a continuous martingale by
Proposition 4.1. All the results of the previous subsection apply
and we are able to compute Vanilla option and CDS prices
conditionally on the knowledge of the law of $H_t$. As a result, we
showed that a CESV model is in fact a timed-changed power of Bessel
process where the subordinator is an integrated time change $H_t=
\int_0^t h_s ds$ with a specific rate of time change $h_t$ that is
defined by
\begin{equation*}
h_t=(1-\alpha)^2\sigma^2_t e^{-2(1-\alpha)r t}
\end{equation*}\\

We now provide two examples of well-known stochastic volatility
models where we compute the law of the time change.

\paragraph{Heston Model} Let us first consider a CIR (1985) diffusion for
the volatility process
\begin{equation*}
d\sigma^2_t=\kappa(\theta-\sigma_t^2)dt+\eta\sigma_t
dW^{\sigma}_t\indent \textmd{and}\indent \sigma^2_0=x>0
\end{equation*}
where $\kappa$, $\theta$ and $\eta$ are strictly positive constants and $W^{\sigma}$ is a
Brownian motion independent from $W$. In fact we are proposing a variation of the Heston (1993) model
by considering $\alpha\neq1$. We then want to compute the law of
\begin{equation*}
H_t=(1-\alpha)^2\int_{0}^{t}\sigma^2_s e^{-2(1-\alpha) r s}ds
\end{equation*}
More precisely, we will compute its Laplace transform, that is to say, for any $\lambda>0$
\begin{equation*}
\mathbb{E}[e^{-\lambda H_t}]
\end{equation*}
For this purpose, let us use the following result:

\newtheorem{laplace1}[BCIR]{Lemma}
\begin{laplace1}
If $X$ a squared Bessel process $BESQ_{(\delta,x)}$ starting from $x\neq0$ and of dimension $\delta$,
then for any function $f:\mathbb{R}_{+}\rightarrow\mathbb{R}_{+}$ such that for any $t>0$: $\int_0^t f(s)ds<\infty$,
we have
\begin{equation*}
\mathbb{E}\bigg[\exp\big(-\int_0^t X_s f(s) ds\big) \bigg]=\frac{1}{\psi'_f(t)^{\delta/2}}\exp\frac{x}{2}\bigg(\phi'_f(0)-\frac{\phi'_f(t)}{\psi'_f(t)}\bigg)
\end{equation*}
where $\phi_f$ is the unique solution of the Sturm-Liouville equation
\begin{equation*}
\phi''_f(s)=2f(s)\phi_f(s)
\end{equation*}
where $s\in[0;\infty[$, $\phi_f(0)=1$, $\phi_f$ is positive and non-increasing and
\begin{equation*}
\psi_f(t)=\phi_f(t)\int_0^t\frac{ds}{\phi^2_f(s)}
\end{equation*}
\end{laplace1}

\begin{proof}
The proof can be found in Pitman and Yor (1982).
\end{proof}\\

\noindent By Lemma 2.1, we can see that
\begin{equation*}
H_t\overset{d}=\big(\frac{2(1-\alpha)}{\eta}\big)^2\int_0^{\frac{\eta^2}{4\kappa}(e^{\kappa t}-1)}X_u \frac{du}{\big(\frac{4\kappa u}{\eta^2}+1\big)^{2[\frac{(1-\alpha)r}{\kappa}+1]}}
\end{equation*}
where $X$ is a $BESQ_{(\frac{4\kappa\theta}{\eta^2},x)}$. Hence for any $\lambda>0$,
\begin{equation*}
\mathbb{E}[e^{-\lambda H_t}]=\mathbb{E}\bigg[\exp\big(-\int_0^{l(t)}
X_s f_{\lambda}(s) ds\big) \bigg]
\end{equation*}
with
\begin{equation}
l(t)=\frac{\eta^2}{4\kappa}(e^{\kappa
t}-1){\label{tchange}}\indent\textmd{and}\indent f_{\lambda}(t)
=\lambda\big(\frac{2(1-\alpha)}{\eta}\big)^2\big(\frac{4\kappa
u}{\eta^2}+1\big)^{-2[\frac{(1-\alpha)r}{\kappa}+1]}
\end{equation}
Defining $a=8((1-\alpha) / \eta)^2$, $b=4\kappa/\eta^2$ and
 $n=-2(\frac{(1-\alpha)r}{\kappa}+1)$ and using Lemma 4.5 we are brought to the resolution of the following ordinary differential equation
\begin{equation*}
\phi''(x)-a\lambda(bx+1)^n\phi(x)=0
\end{equation*}
Then under the boundary conditions, one obtains (see Polyanin and Zaitsev (2003)):
\begin{eqnarray}
\phi_{\lambda}(x)&=&\sqrt{bx+1}~\frac{\frac{\pi}{\sin(\nu\pi)}I_{1/(n+2)}\big(\frac{2\sqrt{a\lambda}}{b(n+2)}(bx+1)^{(n+2)/2}\big)
}{I_{-1/(n+2)}\big(\frac{2\sqrt{a\lambda}}{b(n+2)}\big)}\nonumber\\
&&+\sqrt{bx+1}\frac{K_{1/(n+2)}\big(\frac{2\sqrt{a\lambda}}{b(n+2)}(bx+1)^{(n+2)/2}\big)}{I_{-1/(n+2)}\big(\frac{2\sqrt{a\lambda}}{b(n+2)}\big)}{\label{philam}}\\
\psi_{\lambda}(x)&=&C_1\sqrt{bx+1}~I_{1/(n+2)}\big(\frac{2\sqrt{a\lambda}}{b(n+2)}(bx+1)^{(n+2)/2}\big)\nonumber\\
&&+C_2\sqrt{bx+1}~K_{1/(n+2)}\big(\frac{2\sqrt{a\lambda}}{b(n+2)}(bx+1)^{(n+2)/2}\big){\label{psilam}}
\end{eqnarray}
 where with using the fact that
$I'_{\nu}(x)K_{\nu}(x)-I_{\nu}(x)K'_{\nu}(x)=1/x$ one has
\begin{eqnarray*}
C_1&=&-\frac{b(n+2)}{2a\lambda}K_{1/(n+2)}\big(\frac{2\sqrt{a\lambda}}{b(n+2)}\big)\\
C_2&=&\frac{b(n+2)}{2a\lambda}I_{1/(n+2)}\big(\frac{2\sqrt{a\lambda}}{b(n+2)}\big)
\end{eqnarray*}
We finally obtain the Laplace transform of $H_t$ for any $\lambda>0$
\begin{equation*}
\mathbb{E}[e^{-\lambda
H_t}]=\frac{1}{\psi'_{\lambda}(l(t))^{\delta/2}}\exp\frac{x}{2}\bigg(\phi'_{\lambda}(0)-\frac{\phi'_{\lambda}(l(t))}{\psi'_{\lambda}(l(t))}\bigg)
\end{equation*}
with $\delta=\frac{4\kappa\theta}{\eta^2}$.\\
A simpler example for the forward contract is provided in Atlan and
Leblanc (2005).

\paragraph{Hull and White Model} Let us now consider the Hull and
White (1987) volatility diffusion that is driven by the following
stochastic differential equation:
\begin{equation*}
\frac{d\sigma^2_t}{\sigma^2_t}=\theta dt +\eta dW^{\sigma}_t
\end{equation*}
where $\theta$ and $\eta$ are positive constants and $W^{\sigma}$ is
a Brownian motion independent from $W$. Then $H$ may be computed
and after some simplifications, we obtain:
\begin{equation}
H_t=\frac{4(1-\alpha)\sigma_0^2}{\eta^2}\int_0^\frac{\eta^2t}{4}ds
e^{2(W^{\sigma}_s+\nu s)}
\end{equation}
where $\nu=\frac{2}{\eta^2}(\theta-\frac{\eta^2}{2}-2(1-\alpha)r)$.\\
If we define $A_t^{\nu}=\int_0^t \exp{2(B_s+\nu s)}ds$ where $B$ is
a Brownian motion, we recognize a typical quantity used for the
pricing of Asian options with analytical formulae. Thus, we can
write
\begin{equation*}
H_t=\frac{4(1-\alpha)\sigma_0^2}{\eta^2}A_{\frac{\eta^2t}{4}}^{\nu}
\end{equation*}
and obtain its law using Yor (1992), more precisely we have $\forall
(u,v)\in\mathbb{R}^2_{+}$:
\begin{equation}
f_{|A_t^{\nu}}(u)=\frac{\exp\big(\frac{\pi^2}{2t}-\frac{\nu^2t}{t}-\frac{1}{2u}\big)}{u^2\sqrt{2\pi^3
t}}\int_{-\infty}^{+\infty}dx
e^{x(\nu+1)}e^{-\frac{e^{2x}}{2u}}\psi_{\frac{e^x}{u}}(t)
\end{equation}
where:
\begin{equation}
\psi_r(v)=\int_0^{\infty}dy\exp(-\frac{y^2}{2v})e^{-r\cosh(y)}\sinh(y)\sin(\frac{\pi
y}{v})
\end{equation}

\subsection{Subordinated Bessel Models}

Another way to build stochastic volatility models is to make time
stochastic. Geman, Madan and Yor (2001) recognize that asset prices
may be viewed as Brownian motions subordinated by a random clock.
The random clock may be regarded as a cumulative measure of the
economic activity as said in Clark (1973) and as estimated in
An\'{e} and Geman (2000). The time must be an increasing process,
thus it could either be a L\'{e}vy subordinator or a time integral
of a positive process. In this paper, we only consider the case of a
time integral because we need the continuity of the time change in
order to compute the first-passage time at $0$ to be able to provide
analytical formulas for CDS prices. More generally, for the purpose
of pricing path-dependent options, one needs the continuity of the
time change in order to simulate increments of the time changed
Bessel process. Consequently, we study the case of a time change
$Y_t$ such as
\begin{equation*}
Y_t=\int_0^t y_s ds
\end{equation*}
where the rate of time change $(y_t,t\geq0)$ is a positive
stochastic process.\\
As we have seen in the previous subsection, considering a stochastic
volatility $(\sigma_t, t\geq0)$ in the CEV diffusion is equivalent
to the following rate of time change
\begin{equation*}
y_t=\frac{\sigma^2_t e^{\frac{2rt}{\delta-2}}}{(2-\delta)^2}
\end{equation*}
where $\delta$ is the dimension of the squared Bessel process.
Hence, in order to provide frameworks where one is able to compute
the law of the time change, we are going to go directly through
different modellings of the rate of time change $y_t$.

\paragraph{Integrated CIR Time change} As a first example, let us consider the case where
$y_t$ solves the following diffusion
\begin{equation*}
dy_t=\kappa(\theta-y_t)dt+\eta\sqrt{y_t}dW^Y_t
\end{equation*}
where $W^Y$ is independent from the driving Bessel process. The Laplace transform of $Y_t$ is then defined for any $\lambda>0$ by :
\begin{eqnarray*}
\mathbb{E}[e^{-\lambda Y_t}]&=&e^{\frac{\kappa^2\theta t}{\eta^2}}\frac{\exp\big(-2\lambda y_0/(\kappa+\gamma\coth(\gamma t/2))\big)}
{\big(\cosh(\gamma t/2)+\frac{\kappa}{\gamma}\sinh(\gamma t/2)\big)^{2\kappa\theta/\eta^2}}\\
\gamma&=&\sqrt{\kappa^2+2\eta^2\lambda}
\end{eqnarray*}

\paragraph{Integrated Ornstein-Uhlenbeck Time Change} We now assume the rate of time change to be the solution of the following SDE
\begin{equation*}
dy_t=-\lambda y_t dt +dz_t
\end{equation*}
where $(z_t;t\geq0)$ is a L\'{e}vy subordinator. Let $\psi_z$ denote
the log characteristic function of the subordinator $z_t$,
then
\begin{equation}
\mathbb{E}[e^{iaY_t}]=\exp\big(i a y_0\frac{1-e^{-\lambda
t}}{\lambda}\big)\exp\bigg(\int_0^{a\frac{1-e^{-\lambda
t}}{\lambda}}\frac{\psi_z(x)}{a-\lambda x}dx\bigg)
\end{equation}

Then we can compute the characteristic function of $Y_t$ for
different subordinators and we present here three examples that one
can find in Carr et al. (2003) for which we recall below the characteristic functions:\\
a) For a process with Poisson arrival rate $\nu$ of positive jumps
exponentially distributed with mean $\mu$, we have a L\'{e}vy
density that is
\begin{equation*}
k_z(x)=\frac{\nu}{\mu}e^{-\frac{x}{\mu}}\mathbf{1}_{\{x>0\}}
\end{equation*}
and a log characteristic function
\begin{equation*}
\psi_z(x)=\frac{ix\nu\mu}{1-ix\mu}
\end{equation*}
then we obtain
\begin{equation}
\int\frac{\psi_z(x)}{a-\lambda
x}dx=\log\bigg(\big(x+\frac{i}{\mu}\big)^{\frac{\nu}{\lambda-i\mu
a}}(a-\lambda x)^{\frac{\nu a\mu}{\lambda a\mu+i\lambda}} \bigg)
\end{equation}
b) Let us consider the first time a Brownian motion with drift $\nu$
reaches $1$. It is well known that this passage time follows the
so-called Inverse Gaussian law which L\'{e}vy density and log
characteristic function are respectively
\begin{eqnarray*}
k_z(x)&=&\frac{e^{-\frac{\nu^2 x}{2}}}{\sqrt{2\pi
x^3}}\mathbf{1}_{\{x>0\}}\\
\psi_z(x)&=&\nu-\sqrt{\nu^2-2ix}
\end{eqnarray*}
and we then get
\begin{eqnarray*}
\int\frac{\psi_z(x)}{a-\lambda
x}dx&=&\frac{2\sqrt{\nu^2-2ix}}{\lambda}+\frac{2\sqrt{\nu^2\lambda-2ia}}{\lambda^{3/2}}\textmd{arctanh}\bigg(\sqrt{\frac{\lambda(\nu^2-2ix)}{\nu^2\lambda-2ia}}\bigg)\\
&&-\frac{\nu\log(a-\lambda x)}{\lambda}
\end{eqnarray*}
c) Finally, recall the Stationary Inverse Gaussian case which
L\'{e}vy density and log characteristic function are
\begin{eqnarray*}
k_z(x)&=&\frac{(1+\nu^2x)e^{-\frac{\nu^2 x}{2}}}{2\sqrt{2\pi
x^3}}\mathbf{1}_{\{x>0\}}\\
\psi_z(x)&=&\frac{i u}{\sqrt{\nu^2-2ix}}
\end{eqnarray*}
From these definitions, we obtain
\begin{equation*}
\int\frac{\psi_z(x)}{a-\lambda
x}dx=\frac{\sqrt{\nu^2-2ix}}{\lambda}-\frac{2ia}{\lambda^{3/2}\sqrt{\nu^2\lambda-2ia}}\textmd{arctanh}\bigg(\sqrt{\frac{\lambda(\nu^2-2ix)}{\nu^2\lambda-2ia}}\bigg)
\end{equation*}

\section{Correlation Adjustment}

\subsection{Introducing some Correlation}
We propose a time-changed Bessel process as in the previous
section with some leverage in order to get more independence between
skewness and credit spreads, with respect to which we add a term that contains a
negative correlation (equal to $\rho$) component between the stock return and
the volatility. Hence, let us consider $z_t$ a $\sigma(h_s,s\leq t)$
adapted positive integrable process such as
\begin{equation*}
\frac{e^{\rho z_t}}{\mathbb{E}[e^{\rho z_t}]}
\end{equation*}
is a martingale and a general integrated time change
$H_t=\int_0^t h_s ds$ such as $\mathbb{E}(H_t)<\infty$ then, we can define the stock
price process as follows
\begin{eqnarray*}
S_t&=&e^{rt}BESQ^{2-\delta}_{H_{t\wedge\tau}}\frac{e^{\rho
z_t}}{\mathbb{E}[e^{\rho z_t}]}\\
\tau&=&\inf\{t>0; S_t=0\}
\end{eqnarray*}
where $BESQ$ is a squared Bessel process of dimension $\delta<2$
starting from $S^{1/(2-\delta)}_0$.\\

Let us first show that the process $(e^{-rt}S_t;t\geq0)$ hence
defined is a martingale. We know from Proposition 4.1 that
$BESQ^{2-\delta}_{H_{t\wedge\tau}}$ is a martingale. Now because of
the independence of the processes $z$ and $BESQ$
\begin{equation*}
<BESQ^{2-\delta}_{H_{\cdot\wedge\tau}},\frac{e^{\rho
z_\cdot}}{\mathbb{E}[e^{\rho z_\cdot}]}>_t=0
\end{equation*}
which ensures that $(e^{-rt}S_t;t\geq0)$ is a local martingale. Let
us show that it is actually a true martingale. For this purpose, let
us recall some results:
\newtheorem{defDL}[BCIR]{Definition}
\begin{defDL}
A real valued process $X$ is of class DL if for every $a>0$, the
family of random variables $X_T\mathbf{1}_{\{T<a\}}$ is uniformly
integrable for all stopping times.
\end{defDL}

\noindent We now state the following property:

\newtheorem{LocDL}[BCIR]{Proposition}
\begin{LocDL}
Let $M_t$  be a local martingale such that $\mathbb{E}|M_0|<\infty$
and such that its negative part belongs to class DL. Then its
negative part is a super-martingale. $M_t$ is a martingale if and
only if $\mathbb{E}[M_t]=\mathbb{E}[M_0]$ for all $t>0$.
\end{LocDL}

\begin{proof}
The proof may be found in Elworthy, Li and Yor (1999).
\end{proof}\\

\noindent All the financial assets being positive, one may use a
simpler result than the previous property the proof of which is left
to the reader:

\newtheorem{LocDL2}[BCIR]{Corollary}
\begin{LocDL2}
Let $M_t$  be a positive local martingale such that
$\mathbb{E}|M_0|<\infty$. Then $M_t$ is a supermartingale and it is
a martingale if and only if $\mathbb{E}[M_t]=\mathbb{E}[M_0]$ for
all $t>0$.
\end{LocDL2}

Consequently to prove that the actualized stock price process is a
martingale with regards to the filtration
$\mathcal{F}_t=\mathcal{R}_{H_t}\vee\sigma(h_s;s\leq t)$, we just
need to show that for any $t>0$
\begin{equation*}
\mathbb{E}[e^{-rt}S_t]=S_0
\end{equation*}
which is the case since
\begin{eqnarray*}
\mathbb{E}[e^{-rt}S_t]&=&\mathbb{E}[BESQ^{2-\delta}_{H_{t\wedge\tau}}\frac{e^{\rho
z_t}}{\mathbb{E}[e^{\rho
z_t}]}]=\mathbb{E}\bigg[\mathbb{E}\big[BESQ^{2-\delta}_{H_{t\wedge\tau}}\frac{e^{\rho
z_t}}{\mathbb{E}[e^{\rho z_t}]}\big|\sigma(h_s;s\leq t)\big]\bigg]\\
&=&\mathbb{E}\bigg[\frac{e^{\rho z_t}}{\mathbb{E}[e^{\rho
z_t}]}\mathbb{E}\big[BESQ^{2-\delta}_{H_{t\wedge\tau}}\big|\sigma(h_s;s\leq
t)\big]\bigg]=\mathbb{E}\big[\frac{e^{\rho z_t}}{\mathbb{E}[e^{\rho
z_t}]}S_0\big]=S_0
\end{eqnarray*}

\subsection{Pricing Credit and Equity Derivatives}

The computation of the probability of default is immediate from
Proposition 4.3 because
\begin{equation*}
\tau=\inf\{t\geq0;S_t=0\}=\inf\{t\geq0;BESQ_{H_t}=0\}
\end{equation*}
and then for any $T>0$
\begin{equation*}
\mathbb{P}(\tau\leq T)=\mathbb{E}\bigg[G(1-\frac{\delta}{2},\frac{S_0^{\frac{2}{2-\delta}}}{H_T})\bigg]
\end{equation*}
where G is the complementary Gamma function.\\

Let us compute the European vanilla option prices. For this purpose, we define
$C^{(\rho)}_0(x,y,\delta,K,T;S_0)$
and $P^{(\rho)}_0(x,y,\delta,K,T;S_0)$:
\begin{eqnarray*}
C^{(\rho)}_0(x,y,\delta,K,T;S_0)&=&S_0\frac{e^{\rho
z_T}}{\mathbb{E}[e^{\rho z_T}]}
Q(\frac{(Ke^{-(rT+\rho y)}\mathbb{E}[e^{\rho z_T}])^{\frac{2}{2-\delta}}}{x},4-\delta,\frac{S_0^{\frac{2}{2-\delta}}}{x})\\&&-Ke^{-rT}(1-Q(\frac{S_0^{\frac{2}{2-\delta}}}{x},2-\delta,\frac{(Ke^{-(rT+\rho y)}\mathbb{E}[e^{\rho z_T}])^{\frac{2}{2-\delta}}}{x})) \\
P^{(\rho)}_0(x,y,\delta,K,T;S_0)&=&Ke^{-rT}Q(\frac{S_0^{\frac{2}{2-\delta}}}{x},2-\delta,\frac{(Ke^{-(rT+\rho
y)}\mathbb{E}[e^{\rho
z_T}])^{\frac{2}{2-\delta}}}{x})\\&&-S_0\frac{e^{\rho
z_T}}{\mathbb{E}[e^{\rho z_T}]} (1-Q(\frac{(Ke^{-(rT+\rho
y)}\mathbb{E}[e^{\rho
z_T}])^{\frac{2}{2-\delta}}}{x},4-\delta,\frac{S_0^{\frac{2}{2-\delta}}}{x}))
\end{eqnarray*}
Then, the knowledge of the joint law $\mu_{H_t, z_t}$ for any $t>0$
enables us to compute the option prices as in the previous section:
\begin{eqnarray*}
C_0&=&\int_{\mathbb{R}_+}\int_{\mathbb{R}_+}C^{(\rho)}_0(x,y,\delta,K,T;S_0)\mu_{H_t, z_t}(dx,dy)\\
P_0&=&\int_{\mathbb{R}_+}\int_{\mathbb{R}_+}P^{(\rho)}_0(x,y,\delta,K,T;S_0)\mu_{H_t,
z_t}(dx,dy)
\end{eqnarray*}

\subsection{Examples}

Let us go through most of the time changes presented previously and
see how we can obtain the joint law of the couple $(H_t,z_t)$.

\paragraph{Integrated CIR Time change} Let us consider the following dynamics
\begin{equation*}
dh_t=\kappa(\theta-h_t)dt+\eta\sqrt{h_t}dW^H_t
\end{equation*}
where $W^H$ and $BESQ$ are independent and the stability condition
$\frac{2\kappa\theta}{\eta^2}>1$ is satisfied. Let us take
\begin{equation*}
z_t=h_t+(\kappa-\frac{\rho\eta^2}{2})H_t
\end{equation*}
or equivalently
\begin{equation*}
\rho z_t=\rho(h_0+\kappa\theta
t)+\rho\eta\int_0^t\sqrt{h_s}dW^H_s-\frac{\rho\eta^2}{2}\int_0^t h_s
ds
\end{equation*}
Hence, it is obvious that
\begin{equation*}
\frac{e^{\rho z_t}}{\mathbb{E}[e^{\rho z_t}]}
\end{equation*}
is a local martingale and it is known that it is a martingale as one may check using the Laplace transform below, that
\begin{equation*}
\mathbb{E}[\exp(\rho\eta\int_0^t\sqrt{h_s}dW^H_s-\frac{\rho\eta^2}{2}\int_0^t
h_s ds)]=1
\end{equation*}
In order to compute credit and equity derivatives prices, we then
compute for any positive $\lambda,\mu$ the Laplace transform of
\begin{equation*}
\mathbb{E}[e^{-\lambda H_t-\mu h_t}]
\end{equation*}
It is well known (see Karatzas and Shreve (1991) or Lamberton and Lapeyre (1995)) that
\begin{eqnarray*}
\mathbb{E}[e^{-\lambda H_t-\mu h_t}]&=&\frac{e^{\frac{\kappa^2\theta t}{\eta^2}}}{{\big(\cosh(\gamma t/2)
+\frac{\kappa+\mu\eta^2}{\gamma}\sinh(\gamma t/2)\big)^{2\kappa\theta/\eta^2}}}\exp\big(-h_0 B(t,\lambda,\mu)\big)\\
B(t,\lambda,\mu)&=&\frac{\mu\big(\gamma\cosh(\frac{\gamma
t}{2})-\kappa\sinh(\frac{\gamma
t}{2})\big)+2\lambda\sinh(\frac{\gamma t}{2})}
{\gamma\cosh(\frac{\gamma t}{2})+(\kappa+\lambda\eta^2)\sinh(\frac{\gamma t}{2})}\\
\gamma&=&\sqrt{\kappa^2+2\eta^2\lambda}
\end{eqnarray*}

\paragraph{Heston CESV with correlation} In the same class of models, let us now construct $z$ in terms of the solution
of the following stochastic differential equation
\begin{equation*}
d\sigma^2_t=\kappa(\theta-\sigma^2_t)dt+\eta\sigma_t dW^H_t\indent
\textmd{and} \indent \sigma^2_0=x.
\end{equation*}
First, $h$ is defined by
\begin{equation*}
h_t=\frac{\sigma^2_t e^{2rt/(\delta-2)}}{(2-\delta)^2}
\end{equation*}
Then, following the same method as in the integrated CIR time change case,
we choose $z$ as:
\begin{equation*}
z_t=h_t+(\kappa-\frac{2r}{\delta-2}-\frac{\rho\eta^2}{2})H_t
\end{equation*}
Consequently, $\frac{e^{\rho z_t}}{\mathbb{E}[e^{\rho
z_t}]}$ is a martingale.

Hence, it remains to evaluate for any positive $t$ the Laplace transform
of $(H_t,z_t)$, that is to say for any positive $\lambda,\mu$
\begin{equation*}
\mathbb{E}[e^{-\lambda H_t-\mu h_t}]
\end{equation*}
In order to compute the above quantity, we use the following result
which extends Lemma 4.6 that one can find in Pitman and Yor (1982).

\newtheorem{laplace2}[BCIR]{Lemma}
\begin{laplace2}
If $X$ a squared Bessel process $BESQ_{(\delta,x)}$ starting from
$x\neq0$ and of dimension $\delta$, then for any functions $f$ and
$g$ $:\mathbb{R}_{+}\rightarrow\mathbb{R}_{+}$ such that for any
$t>0$: $\int_0^t f(s)ds<\infty$, we have
\begin{eqnarray*}
\mathbb{E}\bigg[\exp\big(-\int_0^t X_s f(s) ds - g(t)X_t\big)
\bigg]&=&\frac{1}{(\psi'_f(t)+2g(t)\psi_f(t))^{\delta/2}}\times
\\&& \exp\frac{x}{2}\bigg(\phi'_f(0)-\frac{\phi'_f(t)+2g(t)\phi_f(t)}{\psi'_f(t)+2g(t)\psi_f(t)}\bigg)
\end{eqnarray*}
where $\phi_f$ is the unique solution of the Sturm-Liouville
equation
\begin{equation*}
\phi''_f(s)=2f(s)\phi_f(s)
\end{equation*}
where $s\in[0;\infty[$, $\phi_f(0)=1$, $\phi_f$ is positive and
non-increasing and
\begin{equation*}
\psi_f(t)=\phi_f(t)\int_0^t\frac{ds}{\phi^2_f(s)}
\end{equation*}
\end{laplace2}

Taking
\begin{equation*}
g(t)=\frac{\mu}{(2-\delta)^2}\bigg(\frac{\eta^2}{4\kappa
t+\eta^2}\bigg)^{1+\frac{2r}{\kappa(2-\delta)}}
\end{equation*}
in the above Lemma, we obtain that for any positive $\lambda,\mu$
\begin{eqnarray*}
\mathbb{E}[e^{-\lambda H_t-\mu
h_t}]&=&\frac{1}{(\psi'_{\lambda}(l(t))+2\frac{\mu}{(2-\delta)^2}
e^{-(\kappa+\frac{2r}{2-\delta})
t}\psi_{\lambda}(l(t)))^{\delta/2}}\times
\\&&
\exp\frac{x}{2}\bigg(\phi'_{\lambda}(0)-\frac{\phi'_{\lambda}(l(t))+2\frac{\mu}{(2-\delta)^2}
e^{-(\kappa+\frac{2r}{2-\delta})
t}\phi_{\lambda}(l(t))}{\psi'_{\lambda}(l(t))+2\frac{\mu}{(2-\delta)^2}
e^{-(\kappa+\frac{2r}{2-\delta}) t}\psi_{\lambda}(l(t))}\bigg)
\end{eqnarray*}
where noting $\alpha=\frac{\delta-1}{\delta-2}$, the functions
$\phi_{\lambda}$, $\psi_{\lambda}$ and $l$ are defined respectively
in (\ref{tchange}),(\ref{philam}) and (\ref{psilam}).

\paragraph{Integrated Ornstein-Uhlenbeck Time Change} We consider the stochastic time
change $H_t=\int_0^t h_s ds$ and assume that $(h_t;t\geq0)$ is given by
\begin{equation*}
dh_t=-\lambda h_t dt +dz_t
\end{equation*}
where $(z_t;t\geq0)$ is a L\'{e}vy subordinator. Carr et al. (2003)
compute the characteristic function $\Phi(t,a,b)$ of $(H_t, z_t)$
for any $t>0$ and it is given by
\begin{equation}
\mathbb{E}[e^{iaH_t+i b z_t}]=\exp\big(i a h_0\frac{1-e^{-\lambda
t}}{\lambda}\big)\exp\bigg(\int_b^{b+a\frac{1-e^{-\lambda
t}}{\lambda}}\frac{\psi_z(x)}{a+\lambda b-\lambda x}dx\bigg)
\end{equation}
for any $(a,b)\in\mathbb{R}^2_{+}$ where $\psi_z$ is the log
characteristic function of the subordinator. Let us first notice that
\begin{equation*}
\mathbb{E}[e^{\rho z_t}]=\exp(t\psi_z(-i\rho))
\end{equation*}
We quickly recall the computations of $\Phi(t,a,b)$ for different subordinators:\\
a) For a process with Poisson arrival rate $\nu$ of positive jumps
exponentially distributed with mean $\mu$, we obtain
\begin{equation*}
\int\frac{\psi_z(x)}{a+\lambda b-\lambda
x}dx=\log\bigg(\big(x+\frac{i}{\mu}\big)^{\frac{\nu}{\lambda-i\mu
(a+\lambda b)}}((a+\lambda b)-\lambda x)^{\frac{\nu (a+\lambda
b)\mu}{\lambda (a+\lambda b)\mu+i\lambda}} \bigg)
\end{equation*}
b) For an Inverse Gaussian subordinator of parameter $\nu$, we have
\begin{eqnarray*}
\int\frac{\psi_z(x)}{a+\lambda b-\lambda
x}dx&=&\frac{2\sqrt{\nu^2-2ix}}{\lambda}\\&&+\frac{2\sqrt{\nu^2\lambda-2i(a+\lambda
b)}}{\lambda^{3/2}}\textmd{arctanh}\bigg(\sqrt{\frac{\lambda(\nu^2-2ix)}{\nu^2\lambda-2i(a+\lambda
b)}}\bigg)\\&&-\frac{\nu\log((a+\lambda b)-\lambda x)}{\lambda}
\end{eqnarray*}
c) For the Stationary Inverse Gaussian of parameter $\nu$, we write
\begin{eqnarray*}
\int\frac{\psi_z(x)}{a+\lambda b-\lambda
x}dx&=&\frac{\sqrt{\nu^2-2ix}}{\lambda}\\&&-\frac{2i(a+\lambda
b)}{\lambda^{3/2}\sqrt{\nu^2\lambda-2i(a+\lambda
b)}}\textmd{arctanh}\bigg(\sqrt{\frac{\lambda(\nu^2-2ix)}{\nu^2\lambda-2i(a+\lambda
b)}}\bigg)
\end{eqnarray*}



\newpage
\section{Conclusion}
Twelve continuous stochastic stock price models were built in this paper for
equity-credit modelling purposes, all derived from the Constant
Elasticity of Variance model, and as a consequence from Bessel
processes. They all exploit the ability of Bessel processes to be
positive, for those of dimension lower than $2$ to reach $0$ and for
a certain power of a given Bessel process to be a martingale. We
first propose to add a stochastic volatility diffusion to the CEV
model, then more generally to stochastically time change a Bessel
process in order to obtain a stochastic volatility effect, motivated
by known arguments that go back to Clark (1973). Next, in order to
add some correlation between the stock price process and the
stochastic volatility, we extend our framework by multiplying the
Bessel process by exponentials of the volatility and correcting it
by its mean in accordance with arbitrage considerations to obtain
martingale models that are martingales with respect to the joint
filtration of the time-changed Bessel process and the stochastic
time change itself.
Hence, among the different models proposed based on the CEV with
default model, there were first the Constant Elasticity of
Stochastic Variance ones (CESV) taking a Hull and White (1987)
stochastic volatility as well as a Heston (1993) one. We then
proposed integrated time change models, by considering an integrated
CIR time change and an Integrated Ornstein-Uhlenbeck time change
(see Carr et al. (2003)) with different subordinators for the
process driving the diffusion. We finally added correlation between
stock price returns and volatilities to the models presented
previously and provided quasi-analytical formulas for option
and CDS prices for all of them. Let us note that we discussed the true
and local strict martingale properties of CEV processes, that we
 naturally extended to the time change framework.\\
The models presented and discussed in this paper are not specifically
designed to cope just with Equity-Credit frameworks but they also
can be used for instance for FX-rates hybrid modelling by specifying
stochastic interest rates. We can also note that a Poisson jump to
default process can be added to the CEV-like framework in order to
deal with credit spreads for short-term maturities. Campi,
Polbennikov and Sbuelz (2005) and Carr and Linetsky (2005) precisely
considered a CEV model with deterministic volatilities and hazard
rates. The latest paper can easily be generalized to fit in our
time-changed Bessel frameworks. Since our goal was to concentrate on
continuous diffusions, we leave the addition of a jump to default
for further research.

\end{document}